\documentclass[10pt]{article}

\usepackage[utf8]{inputenc}
\usepackage{amsmath}
\usepackage{amsfonts}
\usepackage{amssymb}
\usepackage{eufrak}
\usepackage{mathrsfs}
\usepackage{dsfont}
\usepackage{fancyhdr}
\usepackage{indentfirst}
\usepackage{graphicx}
\usepackage{newlfont}
\usepackage{amssymb}
\usepackage{amsmath}
\usepackage{latexsym}
\usepackage{amsthm}
\usepackage{mathtools}
\usepackage{nicefrac}
\usepackage{epstopdf}
\usepackage{caption}
\usepackage{color}
\usepackage{dsfont}
\usepackage{comment}
\usepackage{hyperref}

\usepackage{esvect}

\usepackage{pdfsync}
\usepackage{amsmath}                                                        
\numberwithin{equation}{section}
\usepackage{graphicx}                                                       
\usepackage{subfigure}
\usepackage{enumitem}                                                    
\usepackage{hyperref}
\usepackage{verbatim}
\usepackage{xcolor}

\usepackage{amssymb}                                                        
\usepackage[mathscr]{eucal}                                             
\usepackage{cancel}                                                             
\usepackage[normalem]{ulem}                                                                 
\usepackage{pstricks}
\usepackage{rotating}
\usepackage{lscape}
\usepackage[paperwidth=8.5in,paperheight=11in,top=1.00in, bottom=1.00in, left=1.00in, right=1.00in]{geometry}
\usepackage{mathtools}                                                      
\mathtoolsset{showonlyrefs=true}                                    
\usepackage{fixltx2e,amsmath}                                           
\MakeRobust{\eqref}

\usepackage{dsfont}

\usepackage{mathdots}
\usepackage{amsthm}                                                             
\allowdisplaybreaks

\theoremstyle{plain}
\newtheorem{theorem}{Theorem}
\numberwithin{theorem}{section}

\newtheorem{lemma}[theorem]{Lemma}                              

\theoremstyle{definition}
\newtheorem{definition}[theorem]{Definition}
\newtheorem{example}[theorem]{Example}
\newtheorem{notation}[theorem]{Notation}
\newtheorem{remark}[theorem]{Remark}
\newtheorem{assumption}[theorem]{Assumption}

\makeatletter
\DeclareRobustCommand{\cev}[1]{%
  {\mathpalette\do@cev{#1}}%
}
\newcommand{\do@cev}[2]{%
  \vbox{\offinterlineskip
    \sbox\z@{$\m@th#1 x$}%
    \ialign{##\cr
      \hidewidth\reflectbox{$\m@th#1\vec{}\mkern4mu$}\hidewidth\cr
      \noalign{\kern-\ht\z@}
      $\m@th#1#2$\cr
    }%
  }%
} \makeatother

\def \y {{\eta}}
\def \s {{\sigma}}

\def \a {{\alpha}}
\def \b {{\beta}}

\def \R {\mathbb{R}}
\def \p {\partial}

\def \PP {\mathsf{P}}

\newcommand\CC{\mathcal{C}}

\newcommand\dd{\mathrm{d}}

\def \rn  {{\mathbb {R}}^{N}}
\def \R  {{\mathbb {R}}}

\def \n {{\nu}}
\def \m {{\mu}}
\def \y {{\eta}}

\def \p {{\partial}}
\def \a {{\alpha}}
\def \O {{\Omega}}

\def \a {{\alpha}}
\def \b {{\beta}}

\def \G {\Ga}

\def \Ga {{\Gamma}}

\def \mm {\mathfrak{m}}
\def \MM {\mathfrak{M}}

\def \It\^o {It\^o }
\def \s {{\sigma}}

\def \R {{\mathbb {R}}}

\def \n {{\nu}}
\def \v {{\nu}}
\def \m {{\mu}}
\def \y {{\eta}}

\def \O {{\Omega}}
\def \phi {{\varphi}}

\def\l {\lambda}

\def \A {\mathcal{A}}

\def \F {\mathcal{F}}

\def \rn {{\mathbb {R}}^{N}}

\def \Ã  {{\`a }}
\def \Ã¨ {{\`e }}
\def \Ã² {{\`o }}
\def \Ã¹ {{\`u }}

\begin{document}

\title{Well-posedness of kinetic McKean-Vlasov equations}

\author{
Andrea Pascucci\thanks{Dipartimento di Matematica, Universit\`a di Bologna, Bologna, Italy.
\textbf{e-mail}: andrea.pascucci@unibo.it} \and Alessio Rondelli\thanks{Dipartimento di
Matematica, Universit\`a di Bologna, Bologna, Italy. \textbf{e-mail}: alessio.rondelli2@unibo.it}}

\date{This version: \today}

\maketitle

\begin{abstract}
We consider the McKean-Vlasov equation
  $$dX_{t}=b(t,X_{t},[X_{t}])dt+\s(t,X_{t},[X_{t}])dW_{t}$$
where $[X_{t}]$ is the law of $X_{t}$. We specifically consider the {\it kinetic case}, where the
equation is degenerate because the dimension of the Brownian motion $W$ is strictly smaller than
that of the solution $X$, as commonly required in classical models of collisional kinetic theory.
Assuming H\"older continuous coefficients and a weak H\"ormander condition, we prove the
well-posedness of the equation. This result advances the existing literature by filling a crucial
gap: it addresses the previously unexplored case where the diffusion coefficient $\sigma$ depends
on the law $[X_t]$. Notably, our proof employs a simplified and direct argument eliminating the
need for PDEs involving derivatives with respect to the measure argument. A critical ingredient is
the sub-Riemannian metric structure induced by the corresponding Fokker-Planck operator.
\end{abstract}


\bigskip\noindent{\bf Acknowledgements}: The authors are members of INDAM-GNAMPA.

\section{Introduction}\label{intro}
Let $Z_{t}=(X_{t},V_{t})$ be a kinetic particle in $\R^{d}\times\R^{d}$, where $X_{t}$ and $V_{t}$
denote, respectively, its position and velocity. Its evolution is typically described by Newton's
law of motion and, in a random environment, by a {McKean-Vlasov (MKV) SDE} of the form
\begin{equation}\label{eL}
  \begin{cases}
    \dd V_{t}=b(t,Z_{t},[Z_{t}])\dd t+
    \s(t,Z_{t},[Z_{t}])\dd W_{t}, \\
    \dd X_{t}=V_{t}\dd t,
  \end{cases}
\end{equation}
where $W$ is a $d$-dimensional Brownian motion and $[Z_{t}]$ denotes the law of $Z_{t}$. Equation
\eqref{eL} is commonly referred to as a {\it second-order} (or {\it inhomogeneous}) system, in
contrast with {\it first-order} (or {\it homogeneous}) systems, where only the velocity variable
is retained. In the latter setting, the position component is ignored and the problem falls within
the classical framework of non-degenerate, uniformly parabolic SDEs, for which a complete and
well-developed theory is available. The present paper is concerned with the {\it kinetic} case,
where the equation is degenerate, since the Brownian motion acts only on the velocity component
and therefore has dimension strictly smaller than that of the full state variable $Z$.

Whenever it exists, a solution to \eqref{eL} is also called a {MKV diffusion}, or a
non-linear Markov process \cite{MR2680971}. In the deterministic case $\s=0$, \eqref{eL} is
related to the classical Vlasov equation \cite{MR3317577}, one of the basic models in plasma
physics and celestial mechanics. The interest in stochastic kinetic particle systems, however,
goes well beyond their intrinsic mathematical relevance and is largely motivated by applications.
In models of collective behaviour, for instance, the stochastic Cucker-Smale system
\cite{MR3847987} plays a prominent role in the study of flocking and self-organization. In
mathematical finance, equations of the form \eqref{eL} arise in the modelling of path-dependent
contingent claims, such as Asian options, and have more recently appeared in the analysis of local
stochastic volatility models; see, for example, \cite{MR2791231,MR4152640}. This framework is also
motivated by recent applications in energy systems: for instance, \cite{PRV26} considers a
stochastic optimal control problem arising from the operation of a large ensemble of residential
storage devices coordinated by an energy aggregator. Further motivations, together with a detailed
discussion of the related literature, are given in Section \ref{secbib}.
{For a survey of MKV equations, with a focus on kinetic models, see
\cite{pagliarani2026mck}.}

{Despite their importance, only partial well-posedness results for the kinetic MKV equation
\eqref{eL} are currently available; see \cite{MR4718402}, \cite{Hao_Rockner_Zhang_2026},
\cite{Pascucci_Rondelli_Veretennikov_2026} and \cite{RusPagl}. In particular, \cite{RusPagl}
considers a constant diffusion coefficient and allows for a drift that may be discontinuous with
respect to the measure argument.}  A major common limitation of these works is that the diffusion
coefficient $\sigma$ is assumed to be independent of the measure argument.


The aim of this paper is precisely to remove this restriction. To the best of our knowledge,
Theorem \ref{t1} below provides the first well-posedness result for kinetic MKV equations with
diffusion coefficient depending on the law, in a genuinely degenerate regime where the drift is
essential for the propagation of noise; see Assumption \ref{H2}. Allowing the diffusion to depend
on the distribution is important not only from a structural viewpoint, but also for applications,
since such a dependence naturally appears in several models. For instance, it is used in
\cite{MR1737547} to describe grazing collisions of particles, and related considerations also
underlie the models studied in \cite{MR3621429}.

{A further motivation comes from the work of Lacker \cite{MR4595391},
who established quantitative propagation-of-chaos estimates for a broad class of mean field
diffusions. As explained in Remarks 2.11 and 4.5 therein, his method also applies to kinetic
systems with constant diffusion. In this setting, our weak well-posedness result for the limiting
MKV equation can be combined with his estimates to obtain quantitative propagation
of chaos, provided the remaining assumptions of \cite{MR4595391} hold. As pointed out in
Remark 2.12 of the same paper, however, the key entropy estimates break down when the diffusion
coefficient involves interactions.}

Let us finally explain the analytical difficulty behind the problem and the strategy of the proof.
In the degenerate kinetic setting, the transition density of the solution does not enjoy the
regularity familiar from the uniformly parabolic case. In particular, one cannot in general expect
Euclidean differentiability in the position variables, not even at first order. This makes it
essential to work in a functional framework adapted to the geometry of the operator, rather than
in standard isotropic spaces. The relevant geometry is the sub-Riemannian one associated with the
kinetic Fokker-Planck operator, and a central ingredient of our analysis is the use of the
corresponding intrinsic regularity theory. This allows us to exploit the smoothing effects
generated by the transport-diffusion structure even in directions where no direct noise is
present.

This intrinsic viewpoint fits naturally with the duality method introduced in \cite{PR24}, which
yields a new and direct proof of well-posedness for MKV equations, different in spirit from the
existing approaches. In the present kinetic framework, the combination of intrinsic regularity and
duality is particularly effective: it avoids the use of PDEs involving derivatives with respect to
the measure variable and, at the same time, is flexible enough to treat diffusion coefficients
with genuine measure dependence.

\subsection{Assumptions and main results}
Let $T>0$ and $\mathcal{P}(\mathbb{R}^N)$ denote the space of probability measures on $\R^{N}$.
For $t\in[0,T]$, we consider the {MKV} system
\begin{equation}\label{e1}
\begin{cases}
  \dd X^{0}_{t}=\left(B^{0} X_{t}+b(t,X_{t},[X_{t}])\right)\dd t+\Sigma(t,X_{t},[X_{t}])\dd W_{t},\\
  \dd X^{1}_{t}=B^{1} X_t \dd t,
\end{cases}
\end{equation}
where $W$ is a $d$-dimensional Brownian motion ($1\le d\le N$) defined on a filtered probability
space $(\O,\F,\PP,\F_{t})$, and $[X_{t}]$ denotes the law of $X_{t}=(X^{0}_{t},X^{1}_{t})$. In
\eqref{e1}, $B^{0}$, $B^{1}$ are constant matrices of dimension $d\times N$ and $(N-d)\times N$
respectively, and the coefficients are measurable functions
\begin{equation}\label{e14}
  b:[0,T]\times\R^{N}\times\mathcal{P}(\mathbb{R}^N)\longrightarrow\R^{d},\qquad
  \Sigma:[0,T]\times\R^{N}\times\mathcal{P}(\mathbb{R}^N)\longrightarrow\R^{d\times d}.
\end{equation}
The prototype of \eqref{e1}, with $N=2d$, is the kinetic particle system \eqref{eL}.

{For a fixed flow of distributions $\m=(\m_t)_{t\in[0,T]}$, we consider
the associated SDE with frozen measure flow}
\begin{equation}\label{e1lin}
\begin{cases}
  \dd X^{0}_{t}=\left(B^{0} X_t+b(t,X_{t},\m_{t})\right)\dd t+\Sigma(t,X_{t},\m_{t})\dd W_{t},\\
  \dd X^{1}_{t}=B^{1} X_t \dd t,\\
  [X_0]=\mu_0
\end{cases}
\end{equation}
Setting
\begin{equation}\label{coeff}
  b^\mu(t,x):=b(t,x,\mu_t),\qquad \Sigma^\mu(t,x):=\Sigma(t,x,\mu_t),\qquad
  \CC^\mu:=\Sigma^{\m}(\Sigma^{\m})^{\ast},\qquad B:=\begin{pmatrix}
   B^{0} \\
   B^{1}
  \end{pmatrix},
\end{equation}
the characteristic operator of \eqref{e1lin} is $\mathcal{A}^\mu_{t,x}+Y$ where
\begin{align}\label{e9}
 \mathcal{A}^\mu_{t,x} &= \frac{1}{2}\sum_{i,j=1}^d
 \CC_{ij}^\mu(t,x)\partial_{x_ix_j}+\sum_{i=1}^d b_i^\mu(t,x)\partial_{x_i}\hspace{30pt}\text{(diffusion)}\\ \label{Y}
 Y&= \partial_t +
 \left\langle B x,\nabla\right\rangle=\partial_t + \sum_{i,j=1}^N B_{ij}x_j\partial_{x_i}.\hspace{40pt}\text{(drift)}
\end{align}
We now introduce the main structural conditions.
\begin{assumption}[\bf Coercitivity in $\R^{d}$]\label{H1} The $d\times d$ matrix
$\CC:=\Sigma\Sigma^{\ast}$ is uniformly positive definite, that is
\begin{equation}\label{e20}
  \langle \CC(t,x,\m)\y,\y\rangle\ge\l|\y|^{2},\qquad (t,x,\m)\in[0,T]\times\R^{N}\times\mathcal{P}(\mathbb{R}^N),\ \y\in\R^{d},
\end{equation}
for some positive constant $\l$.
\end{assumption}
Note that Assumption \ref{e20} is a non-degeneracy condition in $\mathbb{R}^{d}$, while the
system's dimension $N$ can be strictly greater than $d$. 
Propagation of diffusion is ensured by the following hypoellipticity
condition.
\begin{assumption}[\bf Weak H\"ormander's condition]\label{H2}
The vector fields $\partial_{x_1},\cdots,\partial_{x_d}$ and $Y$ satisfy the H\"ormander's
condition
 $$\text{rank Lie}(\partial_{x_1},\dots,\partial_{x_d},Y)=N+1.$$
\end{assumption}
In Section \ref{sec2}, we recall that Assumption \ref{H2} induces a natural \textit{anisotropic}
distance on $\mathbb{R}^{N}$ associated with problem \eqref{e1}. {The coefficients $b$ and $\Sigma$ are assumed to be H\"older continuous
in the spatial variable with respect to this anisotropic distance and Lipschitz continuous
in the measure variable with respect to suitable dual metrics introduced in Section \ref{sec2}.
The precise conditions are stated in Assumption \ref{H3}.}

Our main result is the following
\begin{theorem}[\bf Weak well-posedness]\label{t1}
Under Assumptions \ref{H1},\ref{H2} and \ref{H3} there exists a unique weak solution of \eqref{e1}
with initial distribution $[X_{0}]=\bar{\m}_{0}$. The solution is a strong Markov and Feller
process, and admits a transition density $p$ satisfying the two-sided Gaussian estimate
\begin{equation}\label{ee1}
  C_-\G^{\l_-}(t,x;s,y)\le p(t,x;s,y)\le C_+\G^{\l_+}(t,x;s,y),\qquad 0\le t<s\le T,\ x,y\in\R^{N},
\end{equation}
for some positive constants $C_{\pm},\l_{\pm}$, where $\G^{\l}$ is the Gaussian kernel in
\eqref{Gau}.
\end{theorem}

\begin{remark}\label{r1}
As a direct consequence of Theorem \ref{t1}, if under Assumptions \ref{H1}, \ref{H2} and \ref{H3},
for any fixed flow of distributions $\mu = (\mu_t)_{t\in[0,T]}$ the {frozen SDE} \eqref{e1lin}
admits a unique strong solution (see, for instance, \cite{MR4498506}), then the MKV equation
\eqref{e1} is strongly well-posed.
\end{remark}

\subsection{Motivations and related literature}\label{secbib}
{MKV} equations significantly extend the scope of classical SDEs by incorporating the
effects of interactions between particles or agents subject to the same dynamics. The initial
formulations of MKV equations are attributed to McKean \cite{MR233437}, in the context of
nonlinear PDEs, and independently to Vlasov \cite{Vlasov1968}, in the study of plasma dynamics.
The central idea is that in a multi-agent system of many interacting particles, the effect of the
other particles on any single particle can be averaged out in favor of a mean field; this process
has subsequently been termed propagation of chaos.

MKV equations have garnered significant attention across various fields of applied mathematics,
including mathematical finance, demography, statistical mechanics, and mean-field games (see,
e.g., \cite{MR3752669}). More specifically, MKV equations offer a natural framework for explaining
how aggregative phenomena emerge from microscopic descriptions. Notable examples include their
application in fluid dynamics, biology and population dynamics (see \cite{MR4489768},
\cite{MR1108185}).

A substantial body of literature addresses the well-posedness problem for {\it non-degenerate} MKV
equations (i.e., with $N=d$). Funaki \cite{MR762085} studied the martingale problem under
Lyapunov-type conditions and global Lipschitz assumptions. Subsequently, Shiga and Tanaka
\cite{MR787607} proved that strong well-posedness holds when the diffusion is the identity matrix
and the coefficients are bounded and Lipschitz continuous with respect to the total variation
norm. More recently, Mishura and Veretennikov \cite{MR4421344} established weak existence under
relatively minimal hypotheses and strong well-posedness under more restrictive assumptions using
Girsanov's theorem (which requires independence of the diffusion with respect to the law of the
solution). Independently, Chaudru de Raynal and Frikha \cite{MR4035024}, \cite{MR4377993}, and we
\cite{PR24} proved well-posedness for MKV equations with H\"older coefficients. R\"ockner and
Zhang \cite{MR4255229} (for a measure-independent diffusion coefficient) and Zhao \cite{Zhao} (for
the general case) obtained remarkably comprehensive well-posedness results, even in the presence
of a singular drift term but under the assumptions of uniform non-degeneracy (i.e., with $N=d$).
Finally, quantitative estimates for propagation of chaos are established in \cite{MR4338452}.

As previously mentioned, in the degenerate case, weak existence and strong well-posedness have
been studied in \cite{MR4722362}, \cite{RusPagl}, and \cite{Pascucci_Rondelli_Veretennikov_2026}, even for measurable
coefficients and/or singular drift. However, these results are only valid under the assumption
that {\it the diffusion coefficient is measure-independent}, which excludes the crucial case of
interactions mediated by a collision kernel.

\medskip
The remainder of the paper is structured as follows. In Section \ref{sec2}, we introduce the
sub-Riemannian metric structure associated the the kinetic system. In Section \ref{proof2}, we
prove a crucial duality result, which expresses the push-forward operator in terms of pull-back
operators and allows us to formulate a fixed-point problem in the space of continuous flows of
marginals. Section \ref{proof} is dedicated to the proof of Theorem \ref{t1}.


\section{Anisotropic H\"older spaces}\label{sec2}
As already mentioned, the regularity properties induced by the infinitesimal generator of a
kinetic SDE are drastically different from those in the uniformly non-degenerate case. Within this
framework, a vast literature has been developed by numerous authors, starting from Lanconelli and
Polidoro \cite{LanconelliPolidoro1994}, and Lunardi \cite{Lunardi1997}; here, we recall only a few
fundamental results, such as the intrinsic Taylor formula and the Sobolev embedding theorems,
which were recently established in \cite{MR3429628} and \cite{MR4700191}.

According to \cite{LanconelliPolidoro1994}, Assumption \ref{H2} is equivalent to $B^{1}$ having
the form
\begin{equation}\label{e11}
 B^{1}=\begin{pmatrix} B^{1}_{1} & * & \cdots & * & * \\ 0 & B^{1}_{2} & \cdots & * & *
  \\ \vdots & \vdots & \ddots & \vdots & \vdots
  \\ 0 & 0 & \cdots & B^{1}_{q} & *
  \end{pmatrix}
\end{equation}
where the $*$-blocks are arbitrary and $B^{1}_j$ is a $(d_{j-1}\times d_j)$-matrix of rank $d_j$
with
  $$d = d_0\geq d_1\geq\dots\geq d_q\geq 1,\qquad \sum_{i=0}^q d_i=N.$$
For example, in the Langevin equation \eqref{eL}, we simply have $B^{1}=\begin{pmatrix}
  I_{d} & 0
\end{pmatrix}$ where $I_{d}$ is the $d\times d$ identity matrix.

The structure of $B^{1}$ in \eqref{e11} induces a natural definition of anisotropic norm on
$\R^N$, along with related notions of moments and H\"older continuity.
\begin{definition}[\bf Anisotropic norm]\label{d1}
For any $x\in\R^N$ let
\begin{equation}
|x|_B=\sum_{j=0}^q\sum_{i=\bar{d}_{j-1}+1}^{\bar{d}_j}|x_i|^\frac{1}{2j+1},\qquad
\bar{d}_j=\sum_{k=0}^j d_k,\qquad {\bar{d}_{-1}=0.}
\end{equation}
\end{definition}

\begin{definition}[\bf Anisotropic spaces]
{
Let $\alpha\in]0,1]$:
\begin{itemize}
  \item[i)] $\mathcal{P}^\a_B(\R^N)$ denotes the space of distributions $\m\in\mathcal{P}(\R^N)$ such
  that
\begin{equation}
[\m]_{\mathcal{P}^\a_B}:=\int_{\R^N}|x|^\a_B\m(\dd x)<\infty.
\end{equation}
  \item[ii)] for functions $f:\R^{N}\longrightarrow\R$, we define the seminorm and norm
\begin{equation}
 [f]_{C^\a_B}:=\sup_{x\neq y}\frac{|f(x)-f(y)|}{|x-y|^\a_B},\qquad \|f\|_{bC^\a_B}:= \sup_{x\in\R^{N}}|f(x)|+[f]_{C^\a_B};
\end{equation}
for a function $f:[0,T]\times\R^N\times\mathcal{D}\longrightarrow\R$, where $\mathcal{D}$ is the
domain of the measure variable, we also use the notation
  $$[f]_{C^\a_B}:=\sup_{\substack{t\in[0,T]\\\m\in\mathcal{D}}}[f(t,\cdot,\m)]_{C^\a_B},\qquad
  \|f\|_{bC^\a_B}:=\sup_{\substack{t\in[0,T]\\\m\in\mathcal{D}}}\|f(t,\cdot,\m)\|_{bC^\a_B}.$$
In cases 1) and 2) of Assumption \ref{H3}, we take $\mathcal{D}=\mathcal{P}(\R^N)$ and
$\mathcal{D}=\mathcal{P}^\a_B(\R^N)$, respectively.
\end{itemize}
}
\end{definition}

We now introduce our last main assumption.
\begin{assumption}[\bf Coefficients and initial datum]\label{H3}
The coefficients $b,\Sigma$ in \eqref{e14} are bounded, measurable functions and one of the two
following conditions is satisfied for some $\a>0$:
\begin{itemize}
  \item[1)] the initial law $\bar{\m}_0\in\mathcal{P}(\R^N)$ and there is a positive constant $C$
  such that
\begin{equation}\label{hol1}
 |b(t,x,\m)-b(t,y,\n)| + |\Sigma(t,x,\m)-\Sigma(t,y,\n)|\leq
 C\left(|x-y|^\a_B+\mm^{(1)}_{\a,B}(\m,\n)\right),
\end{equation}
for any $t\in[0,T]$, $x,y\in\R^{N}$ and $\m,\n\in\mathcal{P}(\R^N)$, where
\begin{equation}\label{e25a}
 \mm^{(1)}_{\a,B}(\m,\n):= \sup_{\|f\|_{bC^{\a}_{B}}\leq 1}\int_{\R^{N}}f(x)\left(\mu-\nu\right)(\dd x);
\end{equation}
  \item[2)] the initial law {$\bar{\m}_0$ belongs to the space of distributions with finite Euclidean $\a$-moment,
  i.e. $\bar{\m}_0\in\mathcal{P}_{\text{E}}^\a(\R^N)\subseteq \mathcal{P}_{B}^\a(\R^N)$,}
  and there is a positive constant $C$ such that
\begin{equation}\label{hol2}
  |b(t,x,\m)-b(t,y,\n)| + |\Sigma(t,x,\m)-\Sigma(t,y,\n)|\leq
  C\left(|x-y|^\a_B+\mm^{(2)}_{\a,B}(\m,\n)\right),
\end{equation}
for any $t\in[0,T]$, $x,y\in\R^{N}$ and $\m,\n\in\mathcal{P}^\a_B(\R^N)$, where
\begin{equation}\label{e25b}
 \mm^{(2)}_{\a,B}(\m,\n):= \sup_{[f]_{C^\alpha_B}\leq 1}\int_{\R^{N}} f(x)\left(\mu-\nu\right)(\dd x).
\end{equation}
\end{itemize}
\end{assumption}
\begin{remark}\label{r11}
In the study of the propagation of chaos phenomenon for interacting particle systems (see
\cite{MR1108185}), MKV equations of the form \eqref{e1} emerge, with the coefficients in
\eqref{e14} taking on the normal form
\begin{equation}\label{e13}
  b(t,x,\m)=\int_{\R^{N}}\b(t,x,y)\m(\dd y),\qquad \Sigma(t,x,\m)=\int_{\R^{N}}\s(t,x,y)\m(\dd y).
\end{equation}
{Here $bC^\a_B(\R^N\times\R^N)$ denotes the space of functions $f$
with finite norm
$$
\|f\|_{bC^\a_B(\R^N\times\R^N)}:=
\sup_{x\in\R^N}\|f(x,\cdot)\|_{bC^\a_B}+
\sup_{y\in\R^N}\|f(\cdot,y)\|_{bC^\a_B}.
$$
Condition \eqref{hol1} is satisfied if the kernels $\b$ and $\s$ are measurable and
$$
\sup_{t\in[0,T]}\left(
\|\b(t,\cdot,\cdot)\|_{bC^\a_B(\R^N\times\R^N)}+
\|\s(t,\cdot,\cdot)\|_{bC^\a_B(\R^N\times\R^N)}\right)<\infty.
$$}
\end{remark}
To motivate the use of the metrics $\mm^{(i)}_{\a,B}$, we recall a nice example by Zhao
\cite{Zhao} which shows that having bounded and Lipschitz coefficients with respect to the {\it
total variation distance} is not enough to ensure weak uniqueness even under the hypothesis of
uniformly elliptic diffusion: in fact, a stricter metric is needed. On the other hand, our
approach suggests exploring weaker moduli of continuity than H\"older, a direction left for future
research (see for example \cite{MR4644502} and \cite{MR4466319} where Gaussian estimates have been
proven for Dini continuous coefficients).
%
\begin{example}[\cite{Zhao}, Example 1] Let $W$ be a real Brownian motion, $I$ be the interval $[-2,2]$ and consider the
positive constants
\begin{equation}
  c_1 :=[W_1](I)=\mathcal{N}_{0,1}(I),\qquad c_2:=[2W_1](I)=\mathcal{N}_{0,4}(I),
\end{equation}
and $\l_1,\l_2$ defined as the solutions to the linear system
\begin{equation}
\begin{cases}
  c_1\l_1+(1-c_1)\l_2=1,\\
  c_2\l_1+(1-c_2)\l_2=2.
\end{cases}
\end{equation}
Now, the coefficient
\begin{equation}
 \Sigma(t,\m):=\int_\R\s(t,y)\m(\dd y),\quad \mathrm{where}\quad
 \s(t,y):=\l_1\mathds{1}_{I}\left(\frac{y}{\sqrt{t}}\right)+\l_2\mathds{1}_{I^c}\left(\frac{y}{\sqrt{t}}\right),
\end{equation}
is bounded, uniformly positive and Lipschitz continuous with respect to the total variation
distance. However, both $W_t$ and $2W_t$ are strong solutions to the MKV equation
\begin{equation}
\dd X_t = \Sigma(t,[X_t])\dd W_t,
\end{equation}
with initial value $X_0=0$, since
\begin{equation}
 \Sigma(t,[W_t])=\l_1[W_1](I)+\l_2[W_1](I^c)=1,\qquad
 \Sigma(t,[2W_t])=\l_1[2W_1](I)+\l_2[2W_1](I^c)=2.
\end{equation}
\end{example}

We conclude this section by providing another example with kinetic motivations, where a step-2 Lie
algebra structure arises.
\begin{example} Adding the accelleration component to the Langevin model \eqref{eL}, we get the
following MKV equation in $\R^{3}$ driven by a real Brownian motion $W$: the solution
$Z_{t}=(X_{t},V_{t},A_{t})$ verifies
\begin{equation}\label{eL2}
  \begin{cases}
    \dd A_{t}=
    \s(t,Z_{t},[Z_{t}])\dd W_{t}, \\
    \dd V_{t}=A_{t}\dd t,\\
    \dd X_{t}=V_{t}\dd t,
  \end{cases}
\end{equation}
and the backward Kolmogorov operator
  $$\frac{\s^{2}}{2}\p_{aa}+a\p_{v}+v\p_{x}+\p_{t}$$
satisfies the H\"ormander condition \eqref{H2} because $\p_{a}$ and $Y:=a\p_{v}+v\p_{x}+\p_{t}$,
together with their commutators
  $$[\p_{a},Y]=\p_{v},\qquad [[\p_{a},Y],Y]=\p_{x},$$
span $\R^{4}$. In this case,
  $$B^{1}=\begin{pmatrix}
    1 & 0 & 0 \\
    0 & 1 & 0 \
  \end{pmatrix}$$
induces the anisotropic metric
  $$|(a,v,x)|_{B}=|a|+|v|^{\frac{1}{3}}+|x|^{\frac{1}{5}}.$$
\end{example}

\section{A duality result}\label{proof2}
For a fixed flow of marginals $\m=(\m_{t})_{t\in[0,T]}\in C([0,T];\mathcal{P}(\mathbb{R}^N))$ and
an initial distribution $\bar{\m}_{0}\in \mathcal{P}(\mathbb{R}^N)$, we consider the {frozen SDE}
\eqref{e1lin} with $[X_0]=\bar{\m}_0$ and the corresponding characteristic operator
$\mathcal{A}^\m_{t,x}+Y$ in \eqref{e9}-\eqref{Y}. Note that the integral curve, starting from
$(t,x)$, of the drift $Y$ is the map $s\longrightarrow (s,e^{(s-t)B}x)$.

Since in degenerate problems there are some directions that are more natural than others to
consider and most importantly {\it Euclidean derivatives only exist in some special directions},
we need to properly define the notion of solution to
\begin{equation}\label{e33}
  \mathcal{A}^\m_{t,x} u+Yu=f.
\end{equation}
\begin{definition}[\bf Strong Lie solution]
We say that $u:\,]0,T[\times\R^N\longrightarrow \R$ is a.e. differentiable along $Y$ if there
exists $F\in L^1_{\text{\text{loc}}}(]0,T[,bC(\R^N))$ such that
\begin{equation}\label{e34}
  u(s,e^{(s-t)B}x)=u(t,x)+\int_t^sF(r,e^{(r-t)B}x)\dd r,\qquad (t,x)\in\,]0,T[\times\R^N,\
  s\in\,]0,T[.
\end{equation}
In this case, we set\footnote{Note that, for $F$ continuous, \eqref{e34} amounts to saying that
$F(t,x) = Yu(t,x)$ as a classical directional derivative.} $Yu=F$ and say that $F$ is a Lie
derivative of $u$ along $Y$. 
A solution to equation
  $$\mathcal{A}^\m_{t,x}u+Yu=f$$
on $]0,T[\times\R^N$, is a continuous function $u$ such that $\p_{x_i}u,\p_{x_ix_j}u\in
L^1_{\text{loc}}([0,T[,bC(\R^N))$ for $i,j=1,\dots,d$ and $f-\mathcal{A}^\m_{t,x}u$ is a Lie
derivative of $u$ along $Y$ on $]0,T[\times\R^N$.
\end{definition}
The following result was proven in \cite{MR4660246}.
\begin{theorem}\label{t11}
Under Assumptions \ref{H1},\ref{H2} and \ref{H3}, the operator $\mathcal{A}^\m_{t,x}+Y$ has a
fundamental solution $p^{\m}=p^{\m}(t,x;s,y)$ verifying the following Gaussian estimates
\begin{equation}
  C_-\G^{\l_-}(t,x;s,y)\le p^{\m}(t,x;s,y)\le C_+\G^{\l_+}(t,x;s,y),\qquad 0\le t<s\le T,\
  x,y\in\R^{N},
\end{equation}
where $C_{\pm},\l_{\pm}$ are positive constants which only depend on
$N,d,T,B,\l,\a,\|b\|_{bC^\a_B}$ and $\|\CC\|_{bC^\a_B}$; moreover, $\G^\l$ is the fundamental
solution of $\frac{\l}{2}\Delta+Y$, that is
\begin{equation}\label{Gau}
  \G^\l(t,x;s,y):=G\left(\l\mathbf{C}_{s-t},y-e^{(s-t)B}x\right),
\end{equation}
with
  $$G(\mathbf{C},x):=\frac{1}{\sqrt{(2\pi)^{N}\det\mathbf{C}}}e^{-\frac{1}{2}\langle\mathbf{C}^{-1}x,x\rangle},\qquad
  \mathbf{C}_{s}:=\int_0^s e^{(s-t)B}
  \begin{pmatrix}
    I_d & 0 \\
    0 & 0 \
  \end{pmatrix} e^{(s-t)B^*}\dd t.$$
\end{theorem}

\medskip
{The results in \cite{MR4660246} guarantee that $p^{\m}$ is a strong Lie
solution of the backward Kolmogorov PDE}
\begin{equation}\label{BKPDE}
  (\A^{\m}_{t,x}+Y)p^{\m}(t,x;s,y)=0,\qquad (t,x)\in\,]0,s[\times\R^{N},{}
\end{equation}
{and the frozen SDE \eqref{e1lin} admits a unique weak solution $X^{\m}$
with
transition density $p^{\m}$. 
Applying It\^o's formula to smooth, compactly supported test functions and taking expectations
shows that, for any $(t,x)\in[0,T[\times\R^{N}$, $p^{\m}(t,x;\cdot,\cdot)$ is a distributional
solution of the forward Kolmogorov PDE $(Y^*+\A^{\m,\ast}_{s,y})p^{\m}(t,x;s,y)=0$ on
$]t,T[\,\times\R^{N}$, that is}
\begin{equation}\label{FWPDE}
  \int_{t}^{T}\dd s\int_{\R^{N}}\dd y\,p^{\m}(t,x;s,y)(Y+\A^{\m}_{s,y})\phi(s,y)=0,
\end{equation}
{for any test function $\phi\in C_{0}^{\infty}(]t,T[\,\times\R^{N})$; see
for instance \cite[Section 17.3.1]{PascucciSC}.}

{For simplicity of presentation, in the following computations we
additionally assume continuity in time of the coefficients, so that \eqref{BKPDE} holds pointwise
in the Lie sense\footnote{With merely measurable time dependence,
$p^{\m}(\cdot,\cdot;s,y)$ is absolutely continuous along the integral curves of $Y$. The proof can
then be written using differential identities that hold almost everywhere along these curves.};
see \cite[Remark 1.4]{MR4660246}.}

The {\it push-forward} and the {\it pull-back} operators acting on the distribution
$\y\in\mathcal{P}(\R^{N})$ are defined as
 $$
 \vec{P}^{\m}_{t,s}\y(y)
 :=\int_{\R^{N}}p^\m(t,x;s,y)\y(\dd x),\qquad
 \cev{P}^{\m}_{t,s}\y(x)
 :=\int_{\R^{N}} p^\m(t,x;s,y)\y(\dd y),$$
for $0\le t<s\le T$ and $x,y\in\R^{N}$. Notice that $\vec{P}^{\m}_{0,s}\bar{\m}_{0}$ is the
density of the marginal law $[X^{\m}_{s}]$.
\begin{notation}
Let $\m,\n\in C([0,T];\mathcal{P}(\R^N))$ and $0\le t_1<t_2\le T$. 
For any $\y\in\mathcal{P}(\R^N)$ and $f:\R^N\longrightarrow\R$, suitable measurable function, we
set\footnote{$\vec{I}_{t_1,t_2}^{\m,\n}$ and $\cev{I}_{t_1,t_2}^{\m,\n}$ depend also on
$\bar{\m}_0$ since the push-forward and pull-back operators do. For conciseness, this dependence
is implied throughout the remainder of the paper.}
\begin{align}\label{e23}
 \vec{I}_{t_1,t_2}^{\m,\n}(\y,f)&=\int_{\R^{N}}\dd y\,f(y)(\vec{P}^{\m}_{t_1,t_2}-\vec{P}^{\n}_{t_1,t_2})\y(y),\\ \label{e30}
  \cev{I}_{t_1,t_2}^{\m,\n}(\y,f)&= \int_{\R^{N}}\y(\dd x)\int_{t_1}^{t_2}\dd t\,\cev{P}^{\m}_{t_1,t}
  \left(\A^{\m}_{t,x}-\A^{\n}_{t,x}\right)\cev{P}^{\n}_{t,t_2}f(x).
\end{align}
To fix ideas, $\vec{I}_{t_1,t_2}^{\m,\n}$ and $\cev{I}_{t_1,t_2}^{\m,\n}$ may be interpreted as
forward and backward estimators of the distance between the two flows of distributions $\m$ and
$\n$.
\end{notation}
A crucial tool in our analysis is the following duality result which establishes a link between
the forward and backward estimators $\vec{I}_{t_1,t_2}^{\m,\n}$ and $\cev{I}_{t_1,t_2}^{\m,\n}$.
\begin{lemma}[\bf Duality]\label{l1}
For any $\m,\n\in C([0,T],\mathcal{P}(\R^N))$, $0\le t_1<t_2\le T$, $f\in bC^\a_B(\R^N)$
(re\-specti\-vely, $f\in C^\a_B(\R^{N})$) and $\bar{\m}_0,\y\in\mathcal{P}(\R^N)$ (respectively,
$\bar{\m}_0,\y\in\mathcal{P}^\a_B(\R^N)$), we have
\begin{equation}\label{e100}
  \vec{I}_{t_1,t_2}^{\m,\n}(\y,f)=\cev{I}_{t_1,t_2}^{\m,\n}(\y,f).
\end{equation}
\end{lemma}
\proof We have
\begin{align}
 \vec{I}_{t_1,t_2}^{\m,\n}(\y,f)&=\int_{\R^N}\dd y f(y)\int_{\R^N}\y(\dd x)(p^\m(t_1,x;t_2,y)-p^\n(t_1,x;t_2,y))=
\intertext{(by Fubini's theorem and the Gaussian estimates for $p^{\m}$ and $p^{\n}$, see Theorem
1.1 in \cite{MR4660246})}
  &=\int_{\R^N}\y(\dd x) \int_{\R^N}\dd y f(y)(p^\m(t_1,x;t_2,y)-p^\n(t_1,x;t_2,y))\\
  &=\int_{\R^N}\y(\dd x) \int_{\R^N}\dd y
 f(y)\int_{t_1}^{t_2}\dd t\frac{\dd}{\dd t}\int_{\R^N}\dd
 z (\det e^{tB}) p^\m(t_1,x;t,e^{tB}z)p^\n(t,e^{tB}z;t_2,y)=
\intertext{(by the potential estimates in \cite{MR4660246}, Appendix B, the fact that $\det
e^{tB}=e^{t\,\mathrm{tr}(B)}$, and the definition of Lie derivative along $Y$; we write $Y_{t,z}$
to indicate that $Y$ acts on the variables $(t,z)$)}
  &={\int_{\R^N}\y(\dd x) \int_{\R^N}\dd y f(y)\int_{t_1}^{t_2}\dd t\int_{\R^N}\dd
  z\, \mathrm{tr}(B)e^{t\,\mathrm{tr}(B)}p^\m(t_1,x;t,e^{tB}z)p^\n(t,e^{tB}z;t_2,y)}\\
  &\quad+\int_{\R^N}\y(\dd x) \int_{\R^N}\dd y f(y)\int_{t_1}^{t_2}\dd t\int_{\R^N}\dd
  z \,{{e^{t\,\mathrm{tr}(B)}}\left(Y_{t,z}p^\m(t_1,x;t,e^{tB}z)\right)}p^\n(t,e^{tB}z;t_2,y)\\
  &\quad+\int_{\R^N}\y(\dd x) \int_{\R^N}\dd y f(y)\int_{t_1}^{t_2}\dd t\int_{\R^N}\dd z\,
  \,{e^{t\,\mathrm{tr}(B)}}p^\m(t_1,x;t,e^{tB}z)Y_{t,z}p^\n(t,e^{tB}z;t_2,y)=
\intertext{(by a change of variable)}
 &=  {\int_{\R^N}\y(\dd x) \int_{\R^N}\dd y f(y)\int_{t_1}^{t_2}\dd t\int_{\R^N}\dd
  z\, \mathrm{tr}(B)p^\m(t_1,x;t,z)p^\n(t,z;t_2,y)}\\
  &\quad+\int_{\R^N}\y(\dd x) \int_{\R^N}\dd y f(y)\int_{t_1}^{t_2}\dd t\int_{\R^N}\dd
  z \,\left(Y_{t,z}p^\m(t_1,x;t,z)\right)p^\n(t,z;t_2,y)\\
  &\quad+\int_{\R^N}\y(\dd x) \int_{\R^N}\dd y f(y)\int_{t_1}^{t_2}\dd t\int_{\R^N}\dd z\,
  \,p^\m(t_1,x;t,z)Y_{t,z}p^\n(t,z;t_2,y)=
\intertext{(by Lemma \ref{l3})}
  &= \int_{\R^N}\y(\dd x) \int_{\R^N}\dd y
  f(y)\int_{t_1}^{t_2}\dd t\int_{\R^N}\dd z\, p^\m(t_1,x;t,z)\left(
 \mathcal{A}^\m_{t,z}-\mathcal{A}^\n_{t,z} \right) p^\n(t,z,t_2,y)=\\
\intertext{(using again the potential estimates)}
 &=\int_{\R^N}\y(\dd x) \int_{t_1}^{t_2}\dd t\int_{\R^N}\dd z\, p^\m(t_1,x;t,z)\left(
  \mathcal{A}^\m_{t,z}-\mathcal{A}^\n_{t,z} \right) \int_{\R^N}\dd y f(y) p^\n(t,z,t_2,y)\\
 &=\int_{\R^N}\y(\dd x) \int_{t_1}^{t_2}\dd t\,\cev{P}^{\m}_{t_1,t}\left(\A^{\m}_{t,\cdot}-\A^{\n}_{t,\cdot}\right)\cev{P}^{\n}_{t,t_2}f(x) =\cev{I}^{\m,\n}_{t_1,t_2}(\y,f).\label{e35}
\end{align}
\endproof
\begin{lemma}\label{l3}
Under the same assumptions as in Lemma \ref{l1}, we have that
\begin{equation}\label{e36}
\begin{split}
 &\int_{t_1}^{t_2}\dd t\, \int_{\R^N}\dd z\left[{\left(Y_{t,z}p^\m(t_1,x;t,z)\right)}p^\n(t,z;t_2,y) +
 p^\m(t_1,x;t,z)Y_{t,z}p^\n(t,z;t_2,y)\right]=\\
 &=\int_{t_1}^{t_2}\dd t\, \int_{\R^N}\dd z\, p^\m(t_1,x;t,z)\left( \mathcal{A}^\m_{t,z}-\mathcal{A}^\n_{t,z} \right)
 p^\n(t,z;t_2,y)\\
 &\quad-\mathrm{tr}(B)\int_{t_1}^{t_2}\dd t\, \int_{\R^N}\dd
 z\, p^\m(t_1,x;t,z)p^\n(t,z;t_2,y).
\end{split}
\end{equation}
\end{lemma}
\proof The left-hand side of \eqref{e36} is equal to $J_1+J_2$ where
\begin{align}
 J_{1}&:= \int_{t_1}^{t_2}\dd t\, \int_{\R^N}\dd z \left(Y_{t,z}p^\m(t_1,x;t,z)\right)p^\n(t,z;t_2,y),\\
 J_{2}&:= \int_{t_1}^{t_2}\dd t\, \int_{\R^N}\dd z\, p^\m(t_1,x;t,z)Y_{t,z}p^\n(t,z;t_2,y).
\end{align}
Since $Y=-Y^*-\mathrm{tr}(B)$, we have
\begin{align}
  J_1&= -\int_{t_1}^{t_2}\dd t\, \int_{\R^N}\dd z
  \left[\left(Y_{t,z}^*+\mathrm{tr}(B)\right)p^\m(t_1,x;t,z)\right]p^\n(t,z;t_2,y)\\
  &= -\int_{t_1}^{t_2}\dd t\, \int_{\R^N}\dd z\, p^\m(t_1,x;t,z)Y_{t,z}p^\n(t,z;t_2,y)\\
 &\quad-\mathrm{tr}(B)\int_{t_1}^{t_2}\dd t\, \int_{\R^N}\dd z\,
 p^\m(t_1,x;t,z)p^\n(t,z;t_2,y)=
\intertext{(by the forward Kolmogorov equation)}
  &=\int_{t_1}^{t_2}\dd t\, \int_{\R^N}\dd z\,
 p^\m(t_1,x;t,z)\mathcal{A}^\m_{t,z} p^\n(t,z;t_2,y)\\
 &\quad-\mathrm{tr}(B)\int_{t_1}^{t_2}\dd
 t\, \int_{\R^N}\dd z\, p^\m(t_1,x;t,z)p^\n(t,z;t_2,y).
\end{align}
on the other hand, we have
\begin{align}
 J_2 &= \int_{t_1}^{t_2}\dd t\, \int_{\R^N}\dd z\, p^\m(t_1,x;t,z)Y_{t,z}p^\n(t,z;t_2,y)
\intertext{(by the backward Kolmogorov equation)}
  &=-\int_{t_1}^{t_2}\dd t\, \int_{\R^N}\dd z\, p^\m(t_1,x;t,z)\mathcal{A}^\n_{t,z}p^\n(t,z;t_2,y).
\end{align}
which concludes the proof.
\endproof

\begin{remark}
Lemma \ref{l1} appears to be highly robust to modifications of the problem. For instance, if
$Y=\p_t+\left\langle F(t,x),\nabla \right\rangle$, where $F:[0,T]\times\R^N\rightarrow\R^N$, and
Gaussian as well as potential estimates still hold (as in \cite{MR4498506}), then, due to the
structure of the integral curves, the forward and backward estimators remain equal, and the proof
remains effectively unchanged.
\end{remark}

\section{Proof of Theorem \ref{t1}}\label{proof}
Depending on which of the two alternatives, $i=1$ or $i=2$, in Assumption \ref{H3} is imposed on
the initial law $\bar{\m}_{0}$, the proof of Theorem \ref{t1} relies on the contraction mapping
principle in the metric space
\begin{equation}\label{MM}
  \mathcal{Q}_i:=(C([0,T];\mathcal{P}^\a_B(\R^N)),\MM^{(i)}_{\a,B}),\qquad
  \MM^{(i)}_{\a,B}(\m,\v):=\max\limits_{t\in[0,T]}\mm^{(i)}_{\a,B}(\m_{t},\n_{t}),
\end{equation}
of continuous flows of distributions $\m=(\m_{t})_{t\in[0,T]}$. By adapting the proof of
\cite[Proposition 4.1]{PR24} to the anisotropic metrics $\mm^{(i)}_{\a,B}$, one can readily show
that each $\mathcal{Q}_i$ is complete.

Let $X^{\m}$ denote the solution of \eqref{e1lin} with initial distribution
$[X^{\m}_{0}]=\bar{\m}_{0}$. For each $i=1,2$, under the corresponding alternative in Assumption
\ref{H3}, the map
\begin{align}
  \m\longmapsto([X^{\m}_{t}])_{t\in[0,T]}
\end{align}
maps $\mathcal{Q}_i$ into itself. We further claim that, for sufficiently small $T>0$, this map is
a contraction on $\mathcal{Q}_i$ and therefore admits a unique fixed point $\bar{\m}\in
\mathcal{Q}_i$. This yields the desired conclusion, since $X^{\bar{\m}}$ is then the unique
solution of the MKV equation \eqref{e1}.
%

Recalling the notation \eqref{coeff} for the coefficients of the characteristic operator, by
Assumption \ref{H3} we have
\begin{align}\label{e5}
 \|\CC^{\m}-\CC^{\n}\|_{L^{\infty}([0,T]\times\R^{N})}+\|b^{\m}-b^{\n}\|_{L^{\infty}([0,T]\times\R^{N})}\le
 c\, \MM^{(i)}_{\a,B}(\m,\n),
\end{align}
for some positive constant $c$. Note that
\begin{equation}\label{e22}
  \MM^{(1)}_{\a,B}([X^{\m}],[X^{\n}])=\max_{s\in[0,T]}\sup_{\|f\|_{bC^\alpha_B}\le 1}\vec{I}_{0,s}^{\m,\n}
  (\bar{\m}_0,f),\qquad \MM^{(2)}_{\a,B}([X^{\m}],[X^{\n}])=\max_{s\in[0,T]}\sup_{[f]_{C^\alpha_B}\le 1}\vec{I}_{0,s}^{\m,\n}(\bar{\m}_0,f)
\end{equation}
with $\vec{I}_{0,s}^{\m,\n}(\bar{\m}_0,f)$ as in \eqref{e23}. Now,
$\vec{I}_{0,s}^{\m,\n}(\bar{\m}_0,f)$ can be estimated by means of the duality Lemma \ref{l1}:
indeed, if $[f]_{C^\a_B}\le 1$, we have
\begin{align}
\left|\vec{I}_{0,s}^{\m,\n}(\bar{\m}_0,f)\right|&\stackrel{\ref{l1}}{=}\left|\cev{I}_{0,s}^{\m,\n}(\bar{\m}_0,f)\right|\\
&\le
\int_{0}^{s}\|\cev{P}^{\m}_{0,t}(\A^{\m}_{t,\cdot}-\A^{\n}_{t,\cdot})\cev{P}^{\n}_{t,s}f\|_{L^{\infty}}\dd
t\\ &\le \int_{0}^{s}\|(\A^{\m}_{t,\cdot}-\A^{\n}_{t,\cdot})\cev{P}^{\n}_{t,s}f\|_{L^{\infty}}\dd
t\le \intertext{(by \eqref{e5})} &\le c\, \MM^{(i)}_{\a,B}(\m,\n)\int_{0}^{s}\max_{1\le k,j\le d}
  \left(\|\p_{x_{k}x_{j}}\cev{P}^{\n}_{t,s}f\|_{L^{\infty}}+
  \|\p_{x_{k}}\cev{P}^{\n}_{t,s}f\|_{L^{\infty}}\right)\dd t \le
\intertext{(by the potential estimates for $p^{\n}$, adjusting the constant $c$ as needed, which
depends solely on \newline{$T,\|b\|_{bC^{\a}_{B}}$ and $\|\CC\|_{bC^{\a}_{B}}$, the latter two being
finite by Assumption \ref{H3}})}\label{e7}
  &\le c\, \MM^{(i)}_{\a,B}(\m,\n)\int_{0}^{s}\frac{1}{(s-t)^{1-\frac{\a}{2}}}\dd t\le
  c\, s^{\frac{\a}{2}} \MM^{(i)}_{\a,B}(\m,\n).
\end{align}
Thus, by \eqref{e22}, we have
\begin{equation}
 \MM^{(i)}_{\a,B}([X^{\m}],[X^{\n}])\le cT^\frac{\a}{2}\MM^{(i)}_{\a,B}(\m,\n)
\end{equation}
which shows the contraction property up to small values of $T$. By iterating the argument, we can
extend the result to any finite $T>0$.

To prove the second part of Theorem \ref{t1} we consider a weak solution  $(X,W)$ of \eqref{e1}
defined on a filtered probability space $(\O,\F,\PP,\F_{t})$. The law $\m\equiv \m^{X}$ of $X$ is
uniquely determined, and therefore $X$ solves the {frozen SDE} \eqref{e1lin} with
backward Kolmogorov operator \eqref{BKPDE}; by\footnote{Here we use the fact that the measure
argument enters the infinitesimal generator only through the dependence on $t$, which is merely
measurable, and the results and the bounds in \cite{DiFrancescoPascucci2} are {\it uniform in
$\m$}.} Theorem \ref{t11}, $X$ has transition density $p^{\m}$.

Let $\phi\in bC(\R^{N})$ and
  $$u^{\m}(t,x):=\cev{P}^\m_{t,s}(\phi)=\int_{\R^{N}}p^{\m}(t,x;s,y)\phi(y)dy$$
be the classical solution of the backward Cauchy problem
\begin{equation}\label{ee592}
  \begin{cases}
    (Y+\A_{t}^{\m})u^{\m}(t,x)=0, &\quad (t,x)\in\,]0,s[\times\rn, \\
    u^{\m}(s,x)=\phi(x), &\quad x\in\R^{N}.
  \end{cases}
\end{equation}
The Gaussian estimates for $p^{\m}$ imply that $u^{\m}$ is a bounded function (also uniformly
w.r.t. $\m$) {and has the intrinsic regularity provided by \cite[Theorem
2.7]{luce_pagl_pasc_schauder}. This regularity allows us to apply the intrinsic version of It\^o's
formula, without requiring Euclidean $C^{1,2}$ regularity}. This shows that $u^{\m}(t, X_{t})$ is
a local martingale and a bounded process: therefore, $u^{\m}(t, X_{t})$ is a true martingale and
we have
\begin{equation}\label{ee788}
  \phi(X_{s}) = u^{\m}(t, X_{t}) + \int_{t}^{s} \nabla_{x} u^{\m}(r, X_{r}) \Sigma^{\m}(r, X_{r}) dW_{r}.
\end{equation}
Thus, conditioning \eqref{ee788} on $\mathcal{F}_{t}$, we obtain
  $$\mathbb{E}\left[\phi(X_{s}) \mid \mathcal{F}_{t}\right] = u^{\m}(t, X_{t}) = \int_{\mathbb{R}^{n}} p^{\m}(t, X_{t}; s, y) \phi(y) dy. $$
Given the arbitrariness of $\phi$, it follows that $X$ is a Markov process with transition density
$p^{\m}$: it also follows that $X$ is a Feller process and thus also enjoys the strong Markov
property (see, for instance, Section 18.2 in \cite{PascucciSC}).
\endproof

%
%

\bibliographystyle{acm}
\bibliography{Bibtex-Final}
\end{document}